\documentclass[12pt]{article}
\usepackage{graphicx}

\def\1{{\mathchoice {1\mskip-4mu\mathrm l}
{1\mskip-4mu\mathrm l} {1\mskip-4.5mu\mathrm l}
{1\mskip-5mu\mathrm l}}}
\usepackage{amsmath}
\usepackage{amssymb}
\newcommand{\Z}{\mathbb{Z}}
\newcommand{\R}{\mathbb{R}}
\newcommand{\E}{\mathbb{E}}
\newcommand{\e}{\varepsilon}
\newtheorem{Lemma}{Lemma}[section]
\newtheorem{Definition}{Definition}[section]

\newtheorem{Theorem}{Theorem}[section]
\newtheorem{Example}{Example}[section]
\newtheorem{Corollary}{Corollary}[section]
\newtheorem{itremark}{Remark}[section]

\newenvironment{remark}{\addtocounter{equation}{1}
\begin{itremark}}{\end{itremark}}

\begin{document}

\title{{\bf Short-time Gibbsianness for Infinite-dimensional
Diffusions with Space-Time Interaction}}
\author{
Frank\ Redig
\footnote{Mathematisch Instituut Universiteit Leiden,
Niels Bohrweg 1, 2333 CA Leiden, The Netherlands,
{\em redig@math.leidenuniv.nl}}\\
Sylvie R\oe lly\footnote{Institut f\"ur Mathematik der Universit\"at
Potsdam, Am Neuen Palais 10, 14469 Potsdam, Germany, {\em roelly@math.uni-potsdam.de}}\\
Wioletta Ruszel\footnote{Department of Mathematics and Computing Sciences, University of Groningen, Blauwborgje 3, 9747AC Groningen
{\em W.M.Ruszel@rug.nl}}\\
}
\maketitle

\begin{quote}
Abstract: We consider a class of infinite-dimensional
diffusions where the interaction between the components is both spatial and temporal.
We start the system from a Gibbs measure with finite-range
uniformly bounded interaction. Under suitable conditions
on the drift, we prove that there exists $t_0>0$ such that
the distribution at time
$t\leq t_0$ is a Gibbs measure with absolutely summable interaction.
The main tool is a cluster expansion of both the initial
interaction and certain time-reversed Girsanov factors coming from the
dynamics.
\end{quote}
\vspace{12pt}
\newpage

\section{Introduction}
In this paper we study short-time Gibbsianness for a class of 
infinite-dimensional diffusions with general space-time interaction. 
The diffusion $X=(X_i(t))_{t \geq 0, i \in \Z^d}$ is the solution  of the
following stochastic differential
equations (SDE)
\begin{equation}\label{bobo}
d X_i(t) = b_i(t, X) \, dt + dB_i(t), t > 0, i \in \mathbb{Z}^d
\end{equation}
with values in $\R^{\Z^d}$ and starting at time 0 in a Gibbs measure called $\nu$. 
The drift term $b$ which characterizes the type of interaction between the coordinates 
is adapted but can be non-Markovian, i.e., the interaction at time $t$ may depend on the values of $X$ on the full time interval $[0,t]$.
Moreover, it does not need to be of gradient type.

Earlier short-time Gibbsianness results could be obtained for $b$ defined as a gradient of a Hamiltonian, see
\cite{DerRoe05}. In that case one even has a Gibbs structure on the path space $C(\R_+,\R)^{\Z^d}$ and the existence
of a reversible stationary measure for the dynamics, properties which are heavily used by the authors in their proof. 
In contrast, in our model, 
as soon as the drift is not of gradient form, we do not know if a reversible measure exists
(in fact, even the existence of a stationary measure is not guaranteed).
Another important difference with previous works consists in the fact that, 
since the form of the interaction between the spins (the drift) is quite general, one can not make use any more 
of a decoupling method: this tool was used in \cite{DerRoe05}
to compare the infinite-dimensional dynamics with another simpler dynamics, where a spin is forced to be independent of the others. 
 
For discrete Ising
spins there exist also short-time results; for example, in \cite{efhr}, \cite{LeNRed02}, Gibbsianness is proved for general local
dynamics, making use of a reversible measure. See also \cite{ko}, \cite{er} for results
in this direction for non-discrete bounded spins.
The idea behind short-time conservation of Gibbsianness is that the time-evolved measure
is in a certain sense close to the initial one, which was assumed to be Gibbs. In the case
of discrete spins following a dynamics of Glauber type (see \cite{LeNRed02}), this means that
there is a ``sea'' (in the sense of percolation) of spins that did not
change at all, and isolated islands (for which there is a Peierls estimate)
of spins where at least one flip occurred. This picture is implemented in
a cluster expansion of the Radon-Nikodym derivative 
of the finite-volume distribution at time $t$ w.r.t. 
the finite-volume distribution at time $0$
. In order to obtain Gibbsianness,
one has to show that the sum of the cluster weights containing a fixed
point (the origin) is finite, uniformly in $\Lambda$, for $t$ small enough.
The cluster weights have contributions from the interaction in the
initial measure and from Girsanov factors coming from the dynamics.
The Girsanov factors are multiplicative functionals which are close to 1 
for small $t$, and therefore in good shape for a cluster expansion.

We prove here that for general drifts satisfying assumptions (A1) - (A3) the law of 
the infinite-dimensional diffusion $\eqref{bobo}$ stays Gibbsian for a (short)
time.
In the case of continuous unbounded spins,
the picture is similar as for discrete spins but technically much more envolved.
One has to set up a cluster expansion
of 
the Radon-Nikodym derivative of the finite-volume measure at time $t>0$ too.
The factor coming
from the Girsanov formula now contains stochastic integrals, which
can not be turned into ordinary integrals as is done in the reversible
case, using It\^{o}'s formula. The control of the Girsanov factors reduces to a control of exponential moments
of time-reversal of these stochastic integrals, which can be
done under some regularity conditions.
Furthermore, our results lead, as a corollary, to a constructive 
local existence for a class of infinite-dimensional SDE with non-Markovian drift.

The rest of our paper is organized as follows.
In section \ref{notations} we give definitions of Gibbs measures,
assumptions on the dynamics, as well as some non-trivial examples
satisfying the assumptions (A1) - (A3).
In section \ref{colibri} we
state and prove our main theorem.
In section \ref{raven} we come back to
the examples and verify the assumptions for them.

\section{Notations and Definitions}\label{notations}
In this section we want to define the necessary framework for our study. In particular we present a class of examples
to which our results apply.

\subsection{Interactions and Gibbs measures} \label{interactionsGibbsm}
We will work with measures
on the configuration space $\R^{\Z^d}$. Elements of
$\R^{\Z^d}$ are denoted $x,y,z$. For $\Lambda\subset\Z^d$, $x,y\in \R^{\Z^d}$, we
denote $x_\Lambda z_{\Lambda^c}$ the configuration obtained
by concatenating the restriction of $x$ to $\Lambda $ with
the restriction of $z$ to $\Lambda^c$ ,i,e.,
\begin{equation}
(x_{\Lambda}z_{\Lambda^c})_i = \begin{cases} x_i \text{ if } i \in \Lambda \\ z_i \text{ if } i \in \Lambda^c
\end{cases}.
\end{equation}
We choose the initial distribution $\nu$ to be a \textit{Gibbs measure} associated to an interaction $\varphi$
and an a priori measure $m$. Let us recall some definitions.
\begin{Definition}
An interaction $\varphi$ on $\R^{\Z^d}$ is a collection of functions $\varphi_{\Lambda}$ from $\R^{\Z^d}$ to $\R$, where $\Lambda$ is any finite
subset of $\Z^d$, which satisfy following properties
\begin{enumerate}
 \item $\varphi_{\Lambda}$ is $\mathcal{F}_{\Lambda}$-measurable, where
     $\mathcal{F}_{\Lambda}$ is the sigma-field generated by the canonical projection on $\R^{\Lambda}$.
\item $\varphi$ is absolutely summable, i.e., for all $i\in\Z^d$, 
$$\sum_{\Lambda \ni i}
    ||\varphi_{\Lambda}||_{\infty} < \infty.$$
\item Translation invariance:
\[
 \varphi_{\Lambda +i} (\tau_i x) = \varphi_\Lambda (x)
\]
where $\tau_i$ denotes the shift over $i$: $(\tau_i x)_j:= x_{j-i}$.
\end{enumerate}
\end{Definition}\label{interac}

Furthermore, we assume in this paper that the initial interaction $\varphi$ is 
\begin{itemize}\label{kauwtje}
\item[a)] of finite range, i.e., there exists a $r >0$ such that
if $diam(\Lambda) >r \Rightarrow \varphi_{\Lambda} \equiv 0$ \\
\item[b)] $\forall \Lambda, \varphi_{\Lambda}$ is Lipschitz continuous.
\end{itemize}
Given the interaction $\varphi$ we define the associated \textbf{Hamiltonian} $h=(h_{\Lambda})_{\Lambda \subset \Z^d}$
with respect to the boundary condition $z \in \R^{\Z^d}$ by
\begin{equation} \label{equ:hamil}
 h_{\Lambda}(x_{\Lambda},z_{\Lambda^c})= \underset{\Lambda^{\prime}: \Lambda^{\prime} \cap \Lambda \neq \emptyset}\sum
 \varphi_{\Lambda^{\prime}}(x_{\Lambda}z_{\Lambda^c}).
\end{equation}
(The above sum is finite since $\varphi$ is of finite range. )

The \textbf{finite-volume Gibbs measure} with boundary condition $z$ w.r.t. an a priori measure
 $m$ is
then given by
\begin{equation}
\nu_{\Lambda, z}(dx_{\Lambda}) = \frac{1}{Z_{\Lambda}^z} \exp
(-h_{\Lambda}(x_{\Lambda},z_{\Lambda^c}))\, m(dx_{\Lambda})\label{gibbs}
\end{equation}
where $Z_{\Lambda}^z = \int_{\R^{\Lambda} } \exp (-h_{\Lambda}(y_{\Lambda},z_{\Lambda^c})) m(dy_{\Lambda})$
is the finite-volume partition function. We consider as a priori measure $m$ a finite product measure, absolutely
continuous w.r.t.\ Lebesgue measure.  By absolute summability of the interactions, we then have that the
partition functions
$Z^z_{\Lambda}$ are finite. As usual, the finite-volume Gibbs measure with free boundary condition is defined by
\begin{equation} \label{gibbsfree}
\nu_{\Lambda}(dx_{\Lambda})
:= \frac{1}{Z_{\Lambda}} \exp (-\sum_{A \subset \Lambda} \varphi_A(x_{\Lambda}))\, m(dx_{\Lambda}) 
\end{equation} 
\begin{Definition}
The measure $\mu$ is a Gibbs measure with interaction $\varphi$ and a priori measure $m$, if for all finite $\Lambda \subset
\Z^d$ and smooth test functions $f $ $\mathcal{F}_{\Lambda}$-measurable,
\begin{equation} \label{equ:DLR}
\int f(x_{\Lambda}) \, \mu(dx) = \int \int f(x_{\Lambda}) \, \nu_{\Lambda, z} (dx_{\Lambda}) \mu(dz),
\end{equation}
which means that $\nu_{\Lambda, z} $ is a version of the conditional probability\\
$\mu(dx_{\Lambda}|x_{\Lambda^c}=z_{\Lambda^c})$.
\end{Definition}
Let us recall that 
(\ref{equ:DLR}) is satisfied for all $\Lambda$ as soon as it is satisfied
for singletons (see for example $\cite{Geo88}$), i.e., as soon as for each $i \in \Z^d$ 
\begin{equation}\label{colibrii}
\mu(dx_i|x_{\Z^d \setminus i}) = 
\frac{\exp(-h_i(x_i,x_{\Z^d \setminus i}))}{\int \exp(-h_i(y_i,x_{\Z^d \setminus
i}))m(dy_i)}\,  m(dx_i)
\end{equation}
where $h_i$ is given by (\ref{equ:hamil}) for $\Lambda=\{i\}$.

\begin{remark}
Since we restrict ourselves in this paper to interactions which are uniformly bounded,
certain natural interactions such as quadratic ones are not included here.
In fact, all ``unboundedness'' is hidden in the a priori measure
(the $\log$ of the density of $m$ has to be unbounded for the
partition functions to be finite).
\end{remark}

\subsection{Dynamics}\label{dynamics}

Let $\Omega=C(\R_+, \R)^{\Z^d}$ be the path space for continuous trajectories of the time evolution of the continuous spin system endowed
with the canonical sigma-field $\mathcal{F}$. \\
We denote by $P = \otimes_{i \in \Z^d} P_i$ the Wiener measure on
$\Omega$, resp. by $P^x = \otimes_{i \in \Z^d} P^{x_i}_i$ the Wiener measure with deterministic initial condition
$x=(x_i)_i \in \R^{\Z^d}$, which will be denoted in the finite-volume case just by
$P^x_{\Lambda}=\otimes_{i \in \Lambda}P^{x_i}_i$. \\
Moreover $P^{x,y}_{[0,t]}$ is the law of a Brownian bridge on $[0,t]$
obtained by conditioning $P$ to be at time 0 in $x$ and at time $t$ in $y$. 

The time-reversal $\theta=(\theta_t)_{t > 0}$ is
a family of functionals on the path space $\Omega$. It is defined as follows for $t >  0$ and $\omega(\cdot) \in C([0,t],
\R)^{\Z^d}$:
\begin{equation}
 \theta_t\omega(\cdot ):=\omega(t-\cdot ).
\end{equation}
We consider the following infinite-dimensional system of stochastic differential equations :
\begin{eqnarray}
\begin{cases}
& d X_i(t) = b_i(t,X) \, dt + d B_i(t) , \text{ } t > 0 , i \in \mathbb{Z}^d \\
& X(0) \sim \nu  \label{system1}
\end{cases}
\end{eqnarray}
where $(B_i)_{i \in
\mathbb{Z}^d}$ is a sequence of real-valued independent Brownian motions 
. 
The drift term  $b_i(t,\omega)$ at time $t$ \textit{may possibly}
depend on the values of $\omega$ on the whole time interval $[0,t]$, thus in particular $X$ could be non-Markovian. 
We suppose the existence of a solution
of the system \eqref{system1} and denote it by $Q^\nu$, resp. $Q^{x}$ if the initial condition is deterministic ($\nu=\delta_x$). \\

The drift term $b=(b_i)_i$ satisfies the following assumptions (A1)-(A3).
\begin{itemize}
\item[(A1)] {\bf Translation invariant, finite range and  adapted}, i.e.,
\begin{eqnarray}
\forall i \in \Z^d, b_i(t,\omega) &=& b_0(t, \tau_i \omega) \\
\textrm{and } b_0(t,\omega) = b_0(t,\omega_{\mathcal{N}}) &=& b_0(t,(\omega_{\mathcal{N}}(s) : 0\leq s \leq t ))
\end{eqnarray}
where $\mathcal{N} \subset {\Z^d}$ is a fixed finite connected set containing the origin 
\item[(A2) ] {\bf Lipschitz continuous} uniformly on each compact time interval: \\
$ \forall T>0 \, \exists K(T)>0$ such that
\begin{equation} \label{bLipschitz}
|b_{0}(t,\omega) - b_{0}(t,\omega')| \leq K(T) \sup_{0\leq s \leq t, j \in \mathcal{N}}|\omega_j(s)-\omega'_j(s)|
\end{equation}
\item[(A3) ] {\bf Exponential moment of some time-reversal functional}:\\
$b$ is such that, 
if $F^t_{0}$ denotes the functional 
\begin{equation}
 F^t_{0}(X) := \int_0^t b_0(s,X) \, dX_0(s) - \frac{1}{2} \int_0^t b_0^2(s,X) \, ds
\end{equation} 
its time-reversal $F^t_{0} \circ \theta_t$ is well defined and satisfies 
\begin{equation}
\lim_{t \rightarrow 0} \E_{P^x} \biggl ( \biggl [\exp \bigl (|F^t_{0} \circ
\theta_t |\bigr ) - 1 \biggr ]^{2p} \biggr ) = 0
\end{equation}
for every initial condition $x \in \R^{\Z^d}$.
Here, the integer $p$ is the next odd number larger than $\max(range(b),range(\varphi)) + 1$ 
(see the proof of Lemma \ref{lemmacluster}, equation (\ref{equ:choixdep})).
\end{itemize}

\medskip
We want to present some classes of drifts $b$ which satisfy the above requirements (A1)-(A3).\\

\begin{Example} {\bf Markovian drift}\\
 Let $b_0(t,\omega) = b_0(t,\omega_{\mathcal{N}}(t))$ be a Lipschitz continuous Markovian drift with finite
    range $\mathcal{N}$. Moreover we assume for $j\in \mathcal{N}$ the existence of the first derivative 
$b_0':=\frac{\partial }{\partial x_0}b_0 $ with, for all $T>0$
\begin{equation}
||b||_{T,\infty} + ||b^{\prime}||_{T,\infty} := \underset{t \leq T} \sup \underset{x \in \R^{\Z^d}}\sup
\big(|b_0(t,x)| +  |b'_0(t,x)|\big)< +\infty.
\end{equation}
\end{Example}
This class encloses in particular the Hamiltonian drift from Theorem 1
in \cite{DerRoe05} as a special case. 

\medskip

The second example describes an interaction between the coordinates 
which is spatially degenerate (self-interaction) but has long temporal memory and is thus non-Markovian. 

\begin{Example} {\bf Long memory case}\\
Let the drift $b_i$ be defined as
\begin{equation}
b_i(t,\omega) = \int_0^t \epsilon(s)(\omega_i(s)-\omega_i(0))\, ds
\end{equation}
where the locally integrable \textit{memory function} $\epsilon:[0,\infty] \rightarrow \R $ has the weak continuity property
\begin{equation}\label{contepsilon}
\int_0^t \epsilon(s) ds \underset{t \rightarrow 0} \longrightarrow 0.
\end{equation}
\end{Example}

\medskip

The third example is a generalisation of the first ones.
\begin{Example} {\bf Temporal and spatial interaction}\\
Let $b$ be given by
\begin{equation}
b_i(t,\omega)=\int_0^t \alpha_i(t-s,\omega(s)-\omega(0)) \, dV_s
\end{equation}
where the integrator $V_s$ can be deterministic or stochastic (adapted) and of bounded variation. The functions $\alpha_i$ are
Lipschitz continuous and spatially local, i.e.,
\begin{equation}
\alpha_i(t-s,x) = \alpha_0( t-s,( \tau_ix)_{\mathcal{N}}).
\end{equation}
\end{Example}

The proof that the drifts described in the three examples above satisfy  assumptions (A1)-(A3) is postponed to the Appendix.
In particular we will there explicitly compute the time-reversal functional $F^t_{0} \circ \theta_t$ and provide a proof of the existence of its
exponential moments.


\section{Main result and its proof}\label{colibri}

The following theorem is the main result of our paper.
\begin{Theorem} \label{mainth}
Let us consider an infinite-dimensional Brownian diffusion solution of the system of SDE (\ref{system1})
where the drift term $b$ satisfies the properties (A1)-(A3) given in the previous section. 
Suppose the initial distribution $\nu$ to be a Gibbs measure associated to an a priori measure $m$ and to a finite-range Lipschitz continuous interaction $\varphi$. \\
Then there is a time
$t_0:=t_0(\varphi, b)>0$ such that for any time $t \leq t_0$ the law of the diffusion at time 
 is a Gibbs measure associated to the a priori measure $m$ and to an
absolutely summable interaction $\varphi^t$.
\end{Theorem}
A main step in the proof of Theorem \ref{mainth} is a representation lemma, presented in the next subsection.

\subsection{The finite-dimensional density at time $t$}

Let us first introduce finite-volume dynamics in $\Lambda$ in the following way
\begin{eqnarray} \label{equ:finitevoldyn}
\begin{cases}
& d X^{\Lambda}_i(t) = b_i(t,X) dt + d B_i(t) , \text{ } i \in \Lambda \text{ such that } \mathcal{N} + i \subseteq
\Lambda \\
& d X^{\Lambda}_i(t) = d B_i(t) , \text{ } i \in \Lambda \text{ such that } \mathcal{N} + i \nsubseteq \Lambda.
\end{cases}
\end{eqnarray} 
For the finite-volume Gibbs measure 
$\nu_{\Lambda}$ (the finite-volume measure with free boundary condition) we define 
$\nu^t_{\Lambda}$ to be the distribution of
$(X^{\Lambda}_i(t))_{i \in \Lambda}$ starting from 
$\nu_{\Lambda}$ at time 0.

\begin{Lemma}\label{representlemma}
Let $f^t_{\Lambda}(x_{\Lambda})$ be the density of the finite-volume probability measure $\nu^t_{\Lambda}$ w.r.t. $m(dx_{\Lambda})$. 
 Let therefore $f^0_{\Lambda}$ denote the (Gibbsian) density  of $\nu_\Lambda$ w.r.t. $m(dx_{\Lambda})$.
 Then there exists a time $t_0>0$, such that for any $t \leq t_0$ the ratio
$f^t_{\Lambda}/f^0_{\Lambda}$ admits a cluster representation:
\begin{equation}
\frac{d\nu^t_{\Lambda}}{d\nu_{\Lambda}}(x_{\Lambda}) = \frac{f^t_{\Lambda}(x_{\Lambda})}{f^0_{\Lambda}(x_{\Lambda})} = \exp \biggl ( \sum_{\Gamma \subset \Lambda} a(\Gamma)
w^t(\Gamma, x_{\Gamma}) \biggr ) \label{assump}
\end{equation}
where the``cluster weights" $w^t(\Gamma, x_{\Gamma})$ satisfy
\begin{equation}
\forall i \in \Z^d, \, \sum_{\Gamma \ni i} \sup_{x \in \R^{\Z^d}} |w^t(\Gamma, x_{\Gamma})| < \infty
\end{equation}
and $a(\Gamma)$ are combinatorial factors. 
The sum runs over clusters $\Gamma$ which will be described in $\eqref{clusterw}$.
\end{Lemma}
Suppose Lemma \ref{representlemma} holds true. We now show why it implies the claim of the main theorem when $\nu$ is a Gibbs measure with free boundary counditions. 

Let us denote by $\Upsilon^{t, i}_{\Lambda}(x)$ the conditional density
\begin{equation}
\Upsilon^{t, i}_{\Lambda}(x) := \frac{f^t_{\Lambda}(x_\Lambda )}{\int f^t_{\Lambda}(y_ix_{\Lambda
\setminus i}) m(dy_i)}
\end{equation}
and rewrite it as
\begin{equation} \label{ratio1}
\Upsilon^{t, i}_{\Lambda}(x)=  \frac{f^t_{\Lambda}(x_\Lambda )}{f^0_{\Lambda}(x_\Lambda)} \times \biggl [\int
\frac{f^t_{\Lambda}(y_ix_{\Lambda \setminus i})}{f^0_{\Lambda}(y_ix_{\Lambda \setminus i})}
\frac{f^0_{\Lambda}(y_ix_{\Lambda \setminus i})}{f^0_{\Lambda}(x_ix_{\Lambda \setminus i})} m(dy_i)\biggr
]^{-1}.
\end{equation}
Then using the claim of lemma \ref{representlemma} we have a cluster representation of the ratio
$f^t_{\Lambda}/f^0_{\Lambda}$. So $\eqref{ratio1}$ becomes
\begin{equation}
\begin{split}
& \exp \biggl(\sum_{\Gamma \subset \Lambda} a(\Gamma)w^t(\Gamma, x_{\Gamma}) \biggr) \times \\
& \biggl [ \int \exp \biggl(\sum_{\Gamma \subset \Lambda}a(\Gamma)w^t(\Gamma, (y_ix_{\Lambda \setminus i})_{\Gamma}))
\biggr) \times \exp \biggl( - \sum_{{\Lambda^{\prime} \subset \Lambda}\atop{\Lambda^{\prime} \ni i}}
\varphi_{\Lambda^{\prime}}(y_ix_{\Lambda \setminus i} ) - \varphi_{\Lambda^{\prime}}(x_i x_{\Lambda \setminus i}
)\biggr) m(dy_i)\biggr ]^{-1}\label{clusterw1}
\end{split},
\end{equation}
where we used $\eqref{gibbsfree}$ to express the 
integral in the r.h.s. of
$\eqref{ratio1}$ in terms of the interaction $\varphi$.
The sum runs over all clusters $\Gamma$ whose support is contained in the subset $\Lambda$.
Since the sum over all clusters whose support does not contain $i$ cancel out in the ratio, $\eqref{clusterw1}$
becomes 
\begin{equation}
\Upsilon^{t,i}_{\Lambda}(x) = \frac{\exp \biggl(\underset{{\Gamma \subset \Lambda} \atop {\Gamma \ni i}} \sum
a(\Gamma) w^t(\Gamma, x_{\Gamma}) - \sum_{{\Lambda^{\prime} \subset \Lambda}\atop{\Lambda^{\prime} \ni i}}
\varphi_{\Lambda^{\prime}}(x_\Lambda ) \biggr)}{ \int \exp \biggl(\underset{{\Gamma \subset \Lambda}
\atop {\Gamma \ni i}} \sum a(\Gamma)w^t(\Gamma, (y_ix_{\Lambda \setminus i})_{\Gamma})) -
\sum_{{\Lambda^{\prime}\subset \Lambda}\atop{\Lambda^{\prime} \ni i}} \varphi_{\Lambda^{\prime}}(y_ix_{\Lambda
\setminus i} ) )\biggr) m(dy_i)}.
\end{equation}
Due to the claim of lemma \ref{representlemma}, $\Upsilon^{t,i}_{\Lambda}(x)$ converges, as $\Lambda$ goes to $\Z^d$,
uniformly in $x$, towards
\begin{equation}
\frac{\exp (-h^t_i(x_i,x_{\Z^d \setminus i}))}{\int \exp (-h^t_i(y_i,x_{\Z^d \setminus i})) dm(y_i)}
\end{equation}
where $h^t_i$ is given by
\begin{equation}
h^t_i(x_i,x_{\Z^d \setminus i})= -\sum_{\Gamma \ni i}a(\Gamma) w^t(\Gamma, x_{\Gamma}) + \sum_{\Lambda^{\prime} \ni i}\varphi_{\Lambda^{\prime}}(x)
\end{equation}
which is built from an absolutely summable interaction $\varphi^t$.
In particular we have proven that uniformly in $x$, for $t \leq t_0$
\begin{equation}
d\nu^t_{\Lambda}(x_i|x_{\Lambda \setminus i}) \underset{\Lambda \rightarrow \Z^d} \longrightarrow d\nu^t(x_i|x_{\Z^d
\setminus i})
\end{equation}
which implies that $\nu^t$ is Gibbs, see \eqref{colibrii} in section \ref{interactionsGibbsm}.\\

Let us remark that in order to prove $\eqref{assump}$ in the previous lemma we can replace the
reference measure $m$ by any other one, for example by the Lebesgue measure 
since 
\begin{equation}
\frac{d\nu^t_{\Lambda}}{dm}/\frac{d\nu_{\Lambda}}{dm}=
\frac{d\nu^t_{\Lambda}}{dx}/\frac{d\nu_{\Lambda}}{dx}.
\end{equation}


\subsection{Cluster expansion of the finite-dimensional density}

Let $\nu^t_{\Lambda}$ be the finite-volume time-evolved measure with initial free
boundary condition defined above.
To prove lemma \ref{representlemma} we perform a cluster expansion of $\frac{d\nu^t_{\Lambda}}{dx} /
\frac{d\nu_{\Lambda}}{dx}$.

To do this, we first provide a representation of $f^t_{\Lambda}$ , the density of $\nu^t_{\Lambda}$ 
w.r.t. Lebesgue measure. By the Lebesgue density
theorem, the density $f$ of an absolutely continuous measure $\mu$ w.r.t.\ the Lebesgue measure can be computed via
\begin{equation}
 f(x) = \lim_{\e \rightarrow 0} \int \frac{h_{\e}(y)}{2\e} \mu(dy)
\end{equation}
where $h_{\e}(y)=\1_{[x-\e,x+\e]}(y)$.\\
In a $|\Lambda|$-dimensional situation, one takes $h_{\e}(x_{\Lambda}) = \prod_{i \in \Lambda}h^i_{\e}(x_i)$;\\
Thus, for $\mu=\nu^t_{\Lambda}$
\begin{equation}
  \int \frac{1}{2\e} h_{\e}(x_{\Lambda}) d\nu^t_{\Lambda}(x_{\Lambda}) = 
  \E_{Q_{\Lambda}^{\nu_{\Lambda}}}(\frac{1}{2\e} h_{\e}(X(t)) ) 
=  \int \frac{1}{2\e} \E_{Q_{\Lambda}^{y_{\Lambda}}}(h_{\e}(X(t))
 ) \nu_{\Lambda}(dy_{\Lambda}).\label{densfunct}
\end{equation}
This leads to 
\begin{Lemma}\label{density}
The density $f^t_{\Lambda}$ of $\nu^t_{\Lambda}$ w.r.t. the Lebesgue measure is given by
\begin{eqnarray}
f^t_{\Lambda}(x_{\Lambda})& = & \int_{\R^{\Lambda}} \E_{P^{y,x}_{[0,t],\Lambda}} \biggl(\prod_{i \in \Lambda} \exp
\biggl ( \int_0^t b_i(s,X)d X_i(s) -\frac{1}{2} \int_0^t b_i^2(s,X) ds \biggr ) \times \cr \nonumber\\
&&\prod_{i \in \Lambda}p_t(y_i,x_i) \nu_{\Lambda}(dy_{\Lambda})
\end{eqnarray}
where $p_t$ is the transition kernel of a standard Brownian motion and $P^{y,x}_{[0,t],\Lambda}$ is the law of the
$|\Lambda|$-dimensional Brownian bridge being at time 0 in $y_{\Lambda}$ and in $x_{\Lambda}$ at time $t$.
\end{Lemma}
\textbf{Proof:}\\
Using Girsanov's formula 
\begin{equation}
dQ^{y}_{\Lambda}(X) = \exp \biggl( \sum_{i \in \Lambda} \int_0^t b_i(s,X)dX_i(s) - \frac{1}{2}\int_0^t b^2_i(s,
X)ds \biggr) dP^{y}_{\Lambda}(X)
\end{equation}
 where $P^{y}_{\Lambda}$ is the product of independent Wiener measures. 
Then
\begin{equation}
\begin{split} &\E_{Q_{\Lambda}^{y}}(h_{\e}(X(t)) ) = \cr
& \E_{P_{\Lambda}^{y}} \biggl(\prod_{i \in \Lambda} \exp \biggl ( \int_0^t b_i(s,X)d X_i(s) -\frac{1}{2}
\int_0^t b_i^2(s,X) ds \biggr )h_{\e}^i(X_{i}(t)) \biggr)
\end{split}
\end{equation}
and taking the limit $\e \rightarrow 0$ gives
\begin{equation}
\begin{split}
& \lim_{\e \rightarrow 0} \frac{1}{(2\e)^{|\Lambda|}}\E_{Q_{\Lambda}^{y}}(h_{\e}(X_{\Lambda}(t)) ) = \cr
 & \E_{P_{\Lambda}^{y}} \bigl(\prod_{i \in \Lambda} \exp \bigl ( \int_0^t b_i(s,X)dX_i(s) -\frac{1}{2}
 \int_0^t b_i^2(s,X) ds \bigr ) \bigl \vert X_{\Lambda}(t) = x_{\Lambda} \bigr)\times \prod_{i \in \Lambda}p_t(y_i,x_i) 
\end{split}
\end{equation}
which we can rewrite as
\begin{equation}
 \E_{P^{y,x}_{[0,t],\Lambda}} \bigl(\prod_{i \in \Lambda} \exp \bigl ( \int_0^t b_i(s,X)dX_i(s)
 -\frac{1}{2} \int_0^t b_i^2(s,X) ds \bigr ) \bigr ) \prod_{i \in \Lambda} p_t(y_i,x_i). \label{BB}
\end{equation}
So finally we obtain, plugging in $\eqref{BB}$ into $\eqref{densfunct}$:
\begin{equation}
\begin{split}
f^t_{\Lambda}(x_{\Lambda}) = & \int_{\R^{\Lambda}} \E_{P^{y,x}_{[0,t],\Lambda}} \biggl(\prod_{i \in \Lambda} \exp
\bigl ( \int_0^t b_i(s,X)d X_i(s) -\frac{1}{2} \int_0^t b_i^2(s,X) ds \bigr )\biggr ) \times \cr &
\prod_{i \in \Lambda}p_t(y_i,x_i) \, \nu_{\Lambda}(dy_{\Lambda}).
\end{split}
\end{equation}
\begin{flushright}
 $\square$
\end{flushright}

Using the previous lemma \ref{density} we write the ratio $f^t_{\Lambda}/f^0_{\Lambda}$ as
\begin{equation}\label{ratio2}
\begin{split}
& \frac{f^t_{\Lambda}(x_{\Lambda})}{f^0_{\Lambda}(x_{\Lambda})} = \\
&
\int_{\R^{\Lambda}} \E_{P^{y,x}_{[0,t],\Lambda}} \biggl(\prod_{i \in \Lambda} \exp \biggl ( \int_0^t b_i(s,X)d X_i(s) 
-\frac{1}{2} \int_0^t b_i^2(s,X) ds \biggr ) \times \\
& \prod_{i \in \Lambda} p_t(y_i,x_i)
\exp \biggl(-\sum_{A \subset \Lambda} \varphi_A(y) \biggr)
\exp \biggl(+\sum_{A \subset \Lambda}\varphi_A(x) \biggr) dy_{\Lambda}.
\end{split}
\end{equation}

We will now prove
\begin{Lemma} \label{ratioftf0}
\begin{equation}
\frac{f^t_{\Lambda}(x_{\Lambda})}{f^0_{\Lambda}(x_{\Lambda})} = \E_{P^{x}_{\Lambda}}(R^t_{\Lambda}\circ \theta_t),
\end{equation}
where the functional  $R^t_{\Lambda}$ is defined  as  
\begin{equation}\label{abbrev}
\begin{split}
R^t_{\Lambda}(X):= & \prod_{i \in \Lambda} \exp \biggl (\int_0^t b_i(s,X)dX_i(s) -\frac{1}{2} \int_0^t
b_i^2(s,X)ds\biggr ) \times \\
& \exp \biggl (-\sum_{A \subset \Lambda} \varphi_A(X(0))-\varphi_A(X(t)) \biggr ).
\end{split}
\end{equation}
\end{Lemma}
\textbf{Proof:}\\
Due to $\eqref{abbrev}$ and $\eqref{ratio2}$
\begin{equation}
\frac{f^t_{\Lambda}(x_{\Lambda})}{f^0_{\Lambda}(x_{\Lambda})} = \int_{\R^{\Lambda}} \E_{P^{y,x}_{[0,t],\Lambda}}
\bigl(R^t_{\Lambda} \bigr ) \prod_{i \in \Lambda}p_t(y_i,x_i) \, d y_{\Lambda} .
\end{equation}
Since $\theta_t \circ \theta_t = Id$ we can also write the ratio as
\begin{equation}
\int_{\R^{\Lambda}} \E_{P^{y,x}_{[0,t],\Lambda}} \bigl(R^t_{\Lambda}\circ \theta_t \circ \theta_t \bigr )
\prod_{i \in \Lambda}p_t(y_i,x_i)\, d y_{\Lambda}
\end{equation}
or
\begin{equation}
\int_{\R^{\Lambda}} \E_{P^{y,x}_{[0,t],\Lambda} \circ \theta_t^{-1}} \bigl(R^t_{\Lambda}\circ \theta_t \bigr )
\prod_{i \in \Lambda}p_t(y_i,x_i)\, d y_{\Lambda}.\label{integral}
\end{equation}

The image of the time-reversal of the Brownian bridge is again a Brownian bridge now with reversed starting and final
points, $P^{y,x}_{[0,t],\Lambda} \circ \theta_t^{-1} = P^{x,y}_{[0,t],\Lambda}$.
Furthermore the kernel $p_t$ is symmetric, i.e., $p_t(x,y)=p_t(y,x)$. 
Thus the expectation is now taken w.r.t. a Brownian bridge starting in $x$ and being in $y$ at time $t$.
Integrating out all possible final points $y$, the above integral $\eqref{integral}$ reduces simply to
\begin{equation}
\E_{P^{x}_{\Lambda}}(R^t_{\Lambda}\circ \theta_t) \label{expec}
\end{equation}
which leaves us with an expectation w.r.t. independent Brownian motion starting in $x_{\Lambda}$ of some time-reversed
functional  $R^t_{\Lambda} \circ \theta_t$.
\begin{flushright}
$\square$
\end{flushright}


\subsection{Cluster estimates of $f^t_{\Lambda}/f^0_{\Lambda}$}
To decompose the expectation $\eqref{expec}$ in terms of clusters we write
$R^t_{\Lambda}(X)\circ \theta_t$ under the form $\exp - \sum_{A \subset \Lambda}
\Psi_A(t,X)$ on the path space, where $\Psi_A$ is  $\mathcal{F}_A$-measurable
and apply a standard Mayer expansion for $t$ small.
$\Psi_A$ includes  a contribution from Girsanov terms and the interaction at time 0.

Indeed one can write
\begin{equation}
R^t_{\Lambda}\circ \theta_t (X) = \prod_{A \subset \Lambda} \exp(-\Psi_{A}(t, X))
\end{equation}
with  $\Psi$ is defined as
\begin{equation} \label{eq:PsiA}
 \Psi_{A}(t, X) = \Phi^t_A(X) + \varphi_A(X(t)) - \varphi_A(X(0))
\end{equation}
and
\begin{equation}
 \begin{cases}
   \Phi^t_{\mathcal{N}+i}(X) & = -F^t_{i} \circ \theta_t (X)\quad (F^t_{i}(X):= F^t_{0}(X_{\cdot + i}))\\
\Phi^t_A (X) & \equiv 0 \text{ if there does not exist}\  i\ \text{such that}\  A = \mathcal{N} + i .
 \end{cases}
\end{equation}

Next we give the usual definitions for performing a cluster expansion.
Remember that the drift $b$ and the initial interaction $\varphi$ are of finite range. So we can fix a natural number
$N=N(b,\varphi)$ which depends only on the range of $b$ and $\varphi$ such that for $|A| > N$, $\Psi_{A} \equiv
0$. 
We call a cluster $\gamma = \lbrace A_1,...A_n \rbrace$ a collection of such elements
$A_i$ such that any two $A_i, A_j \in \gamma$ are connected, i.e., there exists a sequence $i=i_1,...,i_m=j$ such
that $A_{i_1} \cap A_{i_2} \neq \emptyset$,...,$A_{i_{m-1}} \cap A_{i_m} \neq \emptyset$. 
The support of the cluster $\gamma$ is the finite subset $\cup_{i=1,\cdots,n} A_i$ and is denoted  by supp($\gamma$). $|\gamma|$
denotes the cardinality of the support of $\gamma$. Clusters $\gamma_i,
\gamma_j$ are said compatible if their supports are disjoint. Let $\mathcal{C}_{\Lambda}$ be the set of all collections of
compatible clusters in $\Lambda$. We expand
\begin{equation}
\prod_{A \subset \Lambda} \bigl( e^{- \Psi_{A}(t,X)} -1 + 1\bigr) = 1+\sum_{n=1}^{\infty}\underset{\lbrace \gamma_1,...\gamma_n
\rbrace \in \mathcal{C}_{\Lambda}} \sum \frac{1}{n!}
\mathcal{K}^{t}(\gamma_1)(X)...\mathcal{K}^{t}(\gamma_n)(X)\label{expand}
\end{equation}
where
\begin{equation}
 \mathcal{K}^{t}(\gamma)(X) = \prod_{A \in \gamma} \biggl ( e^{-\Psi_{A}(t,X)} -1 \biggr ).
\end{equation}
Hence, following Lemma \ref{ratioftf0}, we obtain 
\begin{equation}\label{series}
\begin{split}
f^t_{\Lambda}/f^0_{\Lambda}\, (x_\Lambda)& = \E_{P^{x}_{\Lambda}}\biggl( \prod_{A \subset \Lambda} \bigl(e^{-\Psi_{A}(t,X) } -1 + 1 \bigr) \biggr) \cr
& = 1+ \E_{P^{x}_{\Lambda}}\biggl(\sum_{n=1}^{\infty}\underset{\lbrace \gamma_1,...\gamma_n \rbrace \in
\mathcal{C}_{\Lambda}} \sum \frac{1}{n!}
\mathcal{K}^{t}(\gamma_1)(X)...\mathcal{K}^{t}(\gamma_n)(X)\biggr) \cr
& = 1+ \sum_{n=1}^{\infty}\underset{\lbrace \gamma_1,...\gamma_n \rbrace \in \mathcal{C}_{\Lambda}} \sum \frac{1}{n!}
\E_{P^{x}}\biggl( \mathcal{K}^{t}(\gamma_1)(X)\biggr) ... \E_{P^{x}}\biggl(
\mathcal{K}^{t}(\gamma_n)(X)\biggr) \cr
&=: 1+ \sum_{n=1}^{\infty}\underset{\lbrace \gamma_1,...\gamma_n \rbrace \in \mathcal{C}_{\Lambda}} \sum \frac{1}{n!}
w^t(\gamma_1, x)... w^t(\gamma_n, x)
\end{split}
\end{equation}
where the cluster weights are given by
\begin{equation}
w^t(\gamma, x) := \E_{P^{x}}\biggl(\mathcal{K}^{t}(\gamma)(X)\biggr)
.
\end{equation}
We bound the weights $w^t$ as follows.
\begin{Lemma}\label{lemmacluster}
There exists a strictly positive function $\lambda(t)$ which tends to 0 for $t\to 0$, such that for all
clusters $\gamma$ in $\Lambda$,
\begin{equation}
 \sup_{\Lambda, x }| w^t(\gamma, x) | \leq e^{-c(t)|\gamma|}
\end{equation}
where $c(t):=-\log(\lambda(t))$. 
\end{Lemma}
\textbf{Proof:}\\
The next technical problem is to interchange several times integration and products.
We thus use the following generalised H\"{o}lder inequalities proved in Lemma 5.2 of 
\cite{MinVerZag00}.
\begin{Lemma}\label{MinVer}
Let $(\mu_k)_{k \in \chi}$ be a family of probability measures, each one defined on a space $E_k$ where the indices
$k$ belong to a finite set $\chi$. Let us also define a finite family $(g_i)_i$ of functions on $E_{\chi} = \times_{k
\in \chi} E_k$ such that each $g_i$ is $\chi_i$-local for a certain $\chi_i \subset \chi$ in the sense that
\begin{equation}
g_i(e)=g_i(e_{\chi_i}), \text{ for } e= (e_k)_{k \in \chi} \in E_{\chi}.
\end{equation}
Let $p_i > 1$ be numbers such that 
\[
\forall k \in \chi, \sum_{\{i:\chi_i \ni k\}} 1/p_i \leq 1
\]
Then
\begin{equation}
\biggl | \int_{E_\chi} \prod_i g_i \, \otimes_{k \in \chi} d\mu_k \biggr | \leq \prod_i \biggl ( \int_{E_{\chi_i}} 
|g_i|^{p_i} \otimes_{k \in \chi_i} d\mu_k\biggr )^{1/p_i}
\end{equation}
\end{Lemma}
We apply this lemma with $\chi = $supp$(\gamma)$ ($\gamma=:\lbrace A_1,...,A_n \rbrace$),
$\chi_i=A_i$, $g_i= e^{-\Psi_{A_i}} -1$ and $\mu_k=P^{x}_k$. Let $p > N$ be the next odd number
larger than $N$ and let $p_i=p$ for all $i$. Then $\sum_{A_i \ni k} 1/p_i \leq N/p \leq 1$.
Lemma \ref{MinVer} provides
\begin{equation} 
 | w^t(\gamma, x) | = \biggl | \E_{P^{x}}\biggl(
 \mathcal{K}^{t}(\gamma)(X)\biggr) \biggr | \leq \prod_{i=1}^n \E_{P^{x}_{A_i}}\biggl( \biggl |
 e^{-\Psi_{A_i}(t,X)} -1 \biggr |^p \biggr)^{1/p}. \label{cluster}
\end{equation} \label{equ:choixdep}


Recall that the functional $\Psi_A$ was defined in \eqref{eq:PsiA}.
Since $\varphi_{A}$ is Lipschitz continuous with a constant $C>0$ independent of $A$, the cardinality of $A_i$ is uniformly bounded by $N$ and $\Phi^t_{A} \neq 0$ only if there exists a $k$ such that $A=\mathcal{N} +k$, we obtain
\begin{equation}
|\Psi_{A}(t, X)| \leq   \1_{A=\mathcal{N} +k}|\Phi^t_{\mathcal{N} +k}(X )| + C\sup_{j \in A} |X_j(t)-X_j(0)| . \label{bound1}
\end{equation}
Using the simple fact that, for $a, b \geq 0$,
\begin{equation}
(e^b \cdot e^a - 1)^p \leq 2^p\biggl( e^{p \cdot a}(e^b-1)^p + (e^a-1)^p \biggr )
\end{equation}
 and the estimate \eqref{bound1}, we get
\begin{equation}\label{sum}
\begin{split}
& \E_{P^{x}}\biggl( \bigl |
 e^{-\Psi_{A}(t,X)} -1 \bigr |^p \biggr) \\
& \leq \E_{P^{x}}\biggl( \biggl (\exp \bigl ( |\Phi^t_{\mathcal{N}+k}(X )| \1_{A=\mathcal{N} +k}+ C\sup_{j \in A} |X_j(t)-X_j(0)| \bigr ) -1 \biggr )^p \biggr)  \\
& \leq  2^p \, \E_{P^{x}}\biggl( \exp \bigl (p C\sup_{j \in A} |X_j(t)-X_j(0)| \bigr ) 
\bigl ( \exp (|\Phi^t_{\mathcal{N}+k}(X) |\1_{A=\mathcal{N} +k}) -1 \bigr )^p  \biggr )  \\
& \qquad  + 2^p \E_{P^{x}}\biggl( \bigl( \exp \bigl ( C\sup_{j \in A} |X_j(t)-X_j(0)| \bigr )-1 \bigr)^p \biggl)\\
& =: 2^p \, \E_{P^{x}}(G_{1,A}(t,X)) + 2^p \, \E_{P^{x}}(G_{2,A}(t,X)).
\end{split}
\end{equation}

By the Cauchy-Schwarz inequality 
\begin{equation}
\begin{split}
& \E_{P^{x}}(G_{1,A}(t,X)) \\
& \leq  \E_{P^{x}} \biggl ( \biggl[ \exp \bigr ( \bigl |\Phi^t_{\mathcal{N} + k}(X)|\1_{A=\mathcal{N} +k}\bigr) -1
\biggr ]^{2p} \biggr )^{1/2}
\E_{P^{x}}\biggl( \exp \bigl (2p C \underset{j \in A}\sup |X_j(t)-X_j(0)| \bigr ) \biggr )^{1/2} \\
& \leq  \E_{P^{x}} \biggl ( \bigl( \exp \bigl |F^t_{k}\circ \theta_t (X) | -1 \bigr )^{2p} \biggr )^{1/2}
\E_{P^{x}}\biggl( \exp \bigl (2p C \underset{j \in A}\sup |X_j(t)-X_j(0)| \bigr ) \biggr )^{1/2} \label{factor}
\end{split}
\end{equation}
The exponential moment condition (A4) assures that $ \exp |F^t_{k}(X) \circ \theta_t |$ converges in $\mathbb{L}_{2p}(P^x)$ 
towards 1 for $t$ going to 0 uniformly in $x$.
So there exists a positive function $c_1(t)$ only depending on $t$ vanishing when $t$ is going to 0, such that
\begin{equation}
\E_{P^{x}} \bigl ( \bigl( \exp \bigl |F^t_{k} \circ \theta_t (X)| -1 \bigr )^{2p} \bigr )^{1/2} =: c_1(t).
\end{equation}
The second  term in \eqref{factor}  will be controlled as follows. We recall that $X$ is a family of independent Brownian motions
under $P^x$, thus
\begin{equation} \label{eq:momentexp}
\E_{P^{x}}\biggl( \exp (2p C \underset{j \in A}\sup |X_j(t)-X_j(0)|) \biggr ) \leq \E(\exp(2pCN \sqrt{t}|Z|) ) =: \overline{c}_1(t)
\end{equation}
were $Z$ is a standard Gaussian variable. Clearly, the function $\overline{c}_1(t)$ tends to 1 as $t$ goes to 0. 
We now obtain,
\begin{equation}
\E_{P^{x}}(G_{1,A}(t,X)) \leq c_1(t) \overline{c}_1(t) := C_1(t) \quad \textrm{ with } \lim_{t\rightarrow 0}C_1(t)=0.
\end{equation}
In a similar way we obtain 
\begin{equation}
\begin{split}
 \E_{P^{x}}(G_2(t,p,X)) = & \E_{P^{x}}\biggl( \bigl( \exp \bigl ( C\sup_{j \in A} |X_j(t)-X_j(0)| \bigr )-1 \bigr)^p \biggl)\\
& \leq  \E \biggl( \bigl( \exp \bigl ( C N(b,\varphi)\sqrt{t}|Z| \bigr )-1 \bigr)^p \biggl)\\
& \leq \E \biggl( \bigl ( \int_0^{C N\sqrt{t}|Z|} \exp(u) du \bigr )^p \biggl)\\
 & \leq (C N)^p \sqrt{t}^p \E \biggl ( |Z|^p \exp(pC N\sqrt{t}|Z|)\biggr ) =: c_2(t), 
\end{split}
\end{equation}
where $c_2(t) $ vanishes for $t$ small. 
So finally
\begin{equation}
 \E_{P^{x}}(G_2(t,p,X)) \leq c_2(t) \quad \textrm{ with } \lim_{t\rightarrow 0}c_2(t)=0.
\end{equation}
Thus, calling
\begin{equation}
\lambda(t):= 2\biggl(C_1(t) + c_2(t) \biggr )^{1/Np} \textrm{ and } c(t)= - \log \lambda(t)
\end{equation}
we obtain the desired cluster weight bound
\begin{equation} \label{equ:weightbound}
  | w^t(\gamma,x) | = \bigl | \E_{P^{x}}\bigl(
  \mathcal{K}^t(\gamma)(X)\bigr) \bigr | \leq \exp ( -c(t)|\gamma| ) .
\end{equation}
Note that this bound is uniform in the initial condition $x$.
\begin{flushright}
 $\square$
\end{flushright}
To complete the proof of Lemma \ref{representlemma}, we need a cluster expansion of  $\log(\frac{f^t_{\Lambda}(x_{\Lambda})}{f^0_{\Lambda}(x_{\Lambda})})$.
 This will be done using the Koteck\'{y}-Preiss criterion (see \cite{KotPre86}, p. 492). 
The bound $\eqref{equ:weightbound}$  provides that, for $t \leq t_0$ small enough and any $\gamma \subset \Lambda$,
\begin{equation}
\sup_{x \in \R^{\Z^d}} \sup_{\Lambda \subset \Z^d} \sum_{\gamma^{\prime}: supp(\gamma)\cap supp(\gamma^{\prime})\neq
\emptyset} | w^t(\gamma^{\prime},x) | e^{|\gamma^{\prime}|} \leq |\gamma|.
\end{equation}
Indeed, by the finite-range assumption, the number of clusters $\gamma$ of size $n$, containing
a fixed point is bounded by $e^{cn}$ where $c>0$ does not depend on $t$.
So an absolutely convergent expansion of the
logarithm of the series $\eqref{series}$ exists for $t$ small enough:
\begin{equation}\label{clusterw}
\begin{split}
\log \biggl ( \frac{f^t_{\Lambda}(x_{\Lambda})}{f^0_{\Lambda}(x_{\Lambda})} \biggr ) 
& = \sum_{n=1}^{\infty}\sum_{\Gamma:=\lbrace \gamma_1,...,\gamma_n \rbrace \in U_{\Lambda }}a(\gamma_1,...,\gamma_n)
w^t(\gamma_1,x)... w^t(\gamma_n,x) \\
& =: \sum_{\Gamma \subset \Lambda} a(\Gamma)w^t(\Gamma, x)
\end{split}
\end{equation}
with $a(\Gamma)$ and $a(\gamma_1,...,\gamma_n)$ purely combinatorial terms coming from the Taylor expansion, and $w^t(\Gamma, x)$ depends only on $x_\Gamma$. The set
$U_{\Lambda }$ is the set of all compatible clusters whose union is connected too, the latter sum runs over all
clusters $\Gamma$ which consist of compatible $\gamma_i$.\\
The proof of Lemma \ref{representlemma} is now completed. 
\begin{flushright}
$\square$
\end{flushright}

Next we want to show that if the diffusion starts with any Gibbs measure $\nu$, i.e., not necessarily with a
 measure with free boundary conditions $\nu^{free}$, the probability measure $\nu^t$ is Gibbs associated to the same interaction.
This is done by using the second part of the variational principle characterizing Gibbs measures
(\cite{Geo88}, section 15.4). It applies in our context, even if spins are unbounded,
 since the a priori measure is finite and the interactions are
absolutely summable.

First notice that, if initially the relative entropy density  $i(\nu|\nu^{free})$ vanishes, then the relative entropy density of the time-evolved measure 
satisfies
\[
i(\nu^{t}|\nu^{t,free})\leq i(\nu^1|\nu^{free})=0
\] 
for all $t \geq 0$. Hence if $\nu^{t, free}$ is Gibbs with a absolutely summable interaction, then $\nu^{t}$ is Gibbs with the same interaction. Notice that this fact does not depend on $t$ being small.

In the lemma below we show that $i(\nu|\nu^{free})=0$ is zero for every
extremal Gibbs measure $\nu$ with interaction $\varphi$. By convexity of the relative entropy density, this
then extends to all Gibbs measures $\nu$ with interaction $\varphi$. The proof follows the standard argument of the variational principle
(boundary condition independence of the pressure), see \cite{Geo88}. 
We prefer to spell it out however, for the sake
of completeness, as we are in a context of unbounded spins. 
\begin{Lemma}
Let $\nu$ be an extremal initial Gibbs measure and $\nu^{free}$ a Gibbs measure with free boundary condition. Then
the relative entropy density $i(\nu|\nu^{free})$ is 0.
\end{Lemma}

\textbf{Proof of the lemma:}\\

Let $\Lambda \subset \Z^d$, $\nu$ and $\nu^{free}$ be defined as in the assumption of the lemma. The relative
entropy in volume $\Lambda$ of $\nu$ w.r.t.\ $\nu_{\Lambda}^{free}$ is defined by
\begin{equation}
I_{\Lambda}(\nu_{\Lambda}|\nu_{\Lambda}^{free}) = \int \log \biggl ( \frac{d\nu_{\Lambda}}{d\nu^{free}_{\Lambda}}
(x_{\Lambda}) \biggr ) \nu_{\Lambda}(dx_{\Lambda})
\end{equation}
By the DLR condition,
\begin{equation}
d\nu_{\Lambda}(x_{\Lambda}) = \int \frac{\exp(-h_{\Lambda}(x_{\Lambda},z_{\Lambda^c}))}{Z^z_{\Lambda}}
\nu_{\Lambda}(dz_{\Lambda^c})
\end{equation}
As usual, we show that
\begin{equation}
\frac{d\nu_{\Lambda}}{d\nu^{free}_{\Lambda}}(x_{\Lambda}) \leq \exp (o(|\Lambda|))
\end{equation}
uniformly in $x$, where
\begin{equation}
\frac{d\nu_{\Lambda}}{d\nu^{ free}_{\Lambda}}(x_{\Lambda}) = \int \exp \biggl ( -\sum_{{A \cap \Lambda^c \neq
\emptyset} \atop {A \cap \Lambda \neq \emptyset}} \varphi_A(x_{\Lambda}z_{\Lambda^c})\biggr)
\frac{Z^{free}_{\Lambda}}{Z^z_{\Lambda}}\nu_{\Lambda}(dz_{\Lambda^c}).
\end{equation}
The ratio of the partition functions is equal to
\begin{eqnarray}
\frac{Z^{free}_{\Lambda}}{Z^z_{\Lambda}} & = & \frac{1}{Z^{z}_{\Lambda}} \int \exp( -h_{\Lambda}(x_{\Lambda}))
m(dx_{\Lambda}) \label{ratiopart} \\
& = & \int \exp(-h_{\Lambda}(x_{\Lambda}) + h_{\Lambda}(x_{\Lambda},z_{\Lambda}) )\frac{\exp(
-h_{\Lambda}(x_{\Lambda},z_{\Lambda}))}{Z^{z}_{\Lambda}} m(dx_{\Lambda})\\
& = & \int \exp \biggl ( \sum_{{A \cap \Lambda \neq \emptyset} \atop {A \cap \Lambda^c \neq \emptyset}}
\varphi_{A}(x_{\Lambda}z_{\Lambda^c}) \biggr ) \frac{\exp (- h_{\Lambda}(x_{\Lambda},z_{\Lambda^c})
)}{Z^{z}_{\Lambda}} m(dx_{\Lambda})\\
& = & \E_{\nu^{1}_{\Lambda}} \biggl ( \exp \biggl ( \sum_{{A \cap \Lambda \neq \emptyset} \atop {A \cap \Lambda^c
\neq \emptyset}} \varphi_{A}(X_{\Lambda}z_{\Lambda^c}) \biggr ) \biggr )
\end{eqnarray}
We bound the interaction by its supnorm and use that $\varphi$ is absolute summable to deduce that
\begin{equation}
\sum_{{A \cap \Lambda \neq \emptyset} \atop {A \cap \Lambda^c \neq \emptyset}} \varphi_{A}(x_{\Lambda}z_{\Lambda^c})
\leq o(|\Lambda|)
\end{equation}
which means for the ratio of the partition functions $\eqref{ratiopart}$ that
\begin{equation}
\frac{Z^{free}_{\Lambda}}{Z^z_{\Lambda}} \leq \exp (o(|\Lambda|))
\end{equation}
and a fortiori 
 we conclude that
\begin{eqnarray}
\frac{d\nu_{\Lambda}}{d\nu^{ free}_{\Lambda}}(x_{\Lambda}) & \leq & \int \exp (o(|\Lambda|))
\nu_{\Lambda}(dz_{\Lambda^c}) \\
& = & \exp (o(|\Lambda|)).
\end{eqnarray}
The relative entropy becomes
\begin{equation}
I_{\Lambda}(\nu_{\Lambda}|\nu_{\Lambda}^{free}) \leq \int o(|\Lambda|) \nu_{\Lambda}(dx_{\Lambda}) = o(|\Lambda|)
\end{equation}
and therefore the relative entropy density
\[
i(\nu_{\Lambda}|\nu_{\Lambda}^{free})= \lim_{\Lambda \uparrow\Z^d}\frac{1}{|\Lambda|} I_{\Lambda}(\nu_{\Lambda}|\nu_{\Lambda}^{free})=0
\]
\begin{flushright}
$\square$
\end{flushright}
\begin{remark}
If we want to drop the assumption of translation invariance of the initial interaction $\varphi$ we have to proceed as follows:
First derive as before the cluster
expansion for the free-boundary-condition initial measure; \\then, use the well known fact(see e.g. $\cite{Geo88}$) that every Gibbs measure associated to a given interaction $\varphi$ 
is a mixture of extremal Gibbs measures, which are themselves limits of finite-volume Gibbs measures with fixed boundary conditions.
Now, fix a  boundary condition $z$ and look at the finite-volume dynamics $\eqref{equ:finitevoldyn}$ where the initial distribution is given by  $\nu_{\Lambda, z}$ instead of $\nu_{\Lambda}$. We call $\nu^t_{\Lambda, z}$ the distribution of
$(X^{\Lambda}_i(t))_{i \in \Lambda}$ starting from $\nu_{\Lambda, z}$.
One can without difficulty adapt the result of Lemma \ref{representlemma} to the case with this boundary condition. There exists a similar cluster expansion, with weights (depending on $z$) which can be controlled too. Now the main argument is the following: 
The upper bounds in $\eqref{bound1}$  are uniform in $z$ since the Lipschitz constant $C$ of the interaction is independent of the boundary condition.
Therefore the cluster weights - similar to $\eqref{equ:weightbound}$-  are uniform in $z$, and the cluster expansion -generalizing  $\eqref{clusterw}$- converges when the volume $\Lambda$ goes to $\Z^d$. 
\end{remark}

\begin{Corollary}
The proof of Theorem \ref{mainth} provides a constructive way to obtain a solution of the system $\eqref{bobo}$ on a small time interval as limit (in terms of cluster expansions) of finite-dimensional approximations, whose existence is ensured by the assumption $\eqref{bLipschitz}$. 
\end{Corollary}
\section{Appendix}\label{raven}

In this section we want to show that the assumptions on the drift are satisfied for the presented class of examples.
\subsection{Example 2.1: Markovian Drift} \label{ex:Markov}
We will check the condition (A3) when the drift is Markovian.
First of all we compute the time-reversal of the functional in $X$
\begin{equation}
\int_0^t b_0(s,X(s)) dX_0(s) -  \frac{1}{2} \int_0^t b_0^2(s,X(s)) ds=: I_t(X) - \frac{1}{2} \int_0^t b_0^2(s,X(s)) ds,\label{Markov}
\end{equation}
where the stochastic integral part $I_t(X)$ is defined as
\begin{eqnarray}
I_t(X)=  \int_0^t b_0(s,X(s)) dX_0(s) = \underset{ {n \rightarrow \infty}\atop{ \Delta s \rightarrow 0} } \lim \sum_{j=1}^n
 b_0(s_{j-1},X(s_{j-1}))(X_0(s_j)-X_0(s_{j-1})) \label{driftMarkov}
\end{eqnarray}
with $\Delta s$ the mesh size and $0=s_0 < ...< s_n=t$ a partition of $[0,t]$. 
Then the time-reversal  of the stochastic integral given in 
$\eqref{driftMarkov}$ is
\begin{equation}
\begin{split} 
I_t \circ \theta_t (X) &  = \underset{{n \rightarrow \infty}\atop{ \Delta s \rightarrow 0} } \lim \sum_{j=1}^n
b_0(s_{j-1},X(t-s_{j-1}))(X_0(t-s_j)-X_0(t-s_{j-1})) \\
& \overset{r_j:=t-s_{n-j}} = - \underset{ {n \rightarrow \infty}\atop{ \Delta r \rightarrow 0} } \lim
\sum_{j=1}^n b_0(t-r_{j},X(r_{j}))(X_0(r_j)-X_0(r_{j-1})) \\
& = \underset{ {n \rightarrow \infty}\atop{ \Delta r \rightarrow 0} } \lim \sum_{j=1}^n
b_0(t-r_{j-1},X(r_{j-1}))(X_0(r_j)-X_0(r_{j-1}))  \\
&\quad -  \underset{ {n \rightarrow \infty}\atop{ \Delta r \rightarrow 0} } \lim \sum_{j=1}^n b_0(t-r_{j},X(r_{j}) +
b_0(r_{j-1},X(r_{j-1}))(X_0(r_j)-X_0(r_{j-1})) \\
\end{split}
\end{equation}
which is equal to the sum of an It\^{o} integral and twice a Stratonovich integral,
\begin{equation}
\int_0^t b_0(t-s,X(s)) dX_0(s) - 2\int_0^t b_0(t-s,X(s))\circ dX_0(s).
\end{equation}
Note that $X_0$ is Brownian motion under the measure $P^x$.
So using the It\^{o}-Stratanovich relation (see e.g. definition 3.13 in \cite{KarShr91}), we obtain under $P^x$ 
\begin{equation}
 I_t \circ \theta_t (X)
 = - \int_0^t b_0(t-s,X(s)) dX_0(s) - \int^t_0 b^{\prime}_0(t-s,X(s))ds .
\end{equation}
The second integral in $\eqref{Markov}$ is an ordinary Riemann-Stieltjes integral.
So we obtain
\begin{equation}
\biggl ( \int_0^t b_0^2(s,\cdot(s)) ds \biggr ) \circ \theta_t \, (X) = \int_0^t b_0^2(t-s,X(s)) ds.
\end{equation}

Thus,  the time-reversal of $\eqref{Markov}$ is equal to
\begin{equation}
\begin{split}
& F^t_0\circ \theta_t \,(X) = \\
& - \int_0^t b_0(t-s,X(s)) dX_0(s) - \int^t_0\bigl( b_0^{\prime}(t-s,X(s)) + \frac{1}{2} b_0^2(t-s,X(s))\bigr) ds.
\end{split}
\end{equation}


To obtain the convergence of $\exp |F^t_0\circ \theta_t \,(X)|$ towards 1 when $t$ tends to 0 in $L^{2p}(P^x)$, since the a.s. convergence is clear, it is enough to prove a uniform bound for $t \in [0,1]$ in $L^{2p'}, p'>p$. Indeed
\begin{equation}
\begin{split}
& \E_{P^x}\biggl ( \exp \biggl ( 2p' \text{ }  |F^t_0\circ \theta_t(X) | \biggr ) \biggr )
\\
& \leq e^{p't(||b_0||^2_{\infty} + 2||b_0^{\prime}||_{\infty})}\E_{P^x}\biggl ( \exp\bigl( 2p'  |\int_0^t b_0(t-s,X(s)) dX_0(s)|\bigr)
\biggr ) .
\label{bound2}
\end{split}
\end{equation}
The first term on the right side is bounded for $t \in [0,1]$. The second term can be controlled as follows:
\begin{equation}
\begin{split}
&\E_{P^x}\biggl ( \exp\bigl( 2p'  |\int_0^t b_0(t-s,X(s)) dX_0(s)|\bigr) \biggr ) \\
& \leq \E_{P^x}\biggl ( \exp\bigl( 2p'  \int_0^t b_0(t-s,X(s)) dX_0(s)\bigr) \biggr )
+ \E_{P^x}\biggl ( \exp\bigl( -2p'  \int_0^t b_0(t-s,X(s)) dX_0(s)\bigr) \biggr ).
\label{bound2bis}
\end{split}
\end{equation}
Since $\exp\bigl( 2p'  \int_0^t b_0(t-s,X(s)) dX_0(s) -2 p'^2 \int_0^t b_0^2(t-s,X(s))ds\bigr )$ 
(resp. \\
$\exp\bigl( -2p'  \int_0^t b_0(t-s,X(s)) dX_0(s) -2 p'^2 \int_0^t b_0^2(t-s,X(s))ds\bigr )$)
is a $P^x$-martingale with expectation 1 
\begin{equation}
\begin{split}
\E_{P^x}\biggl ( \exp\bigl( 2p' \int_0^t b_0(t-s,X(s)) dX_0(s)\bigr) \biggr ) & \leq e^{2p'^2t ||b_0||^2_{\infty}}\\
 \textrm{and } \E_{P^x}\biggl ( \exp\bigl( -2p' \int_0^t b_0(t-s,X(s)) dX_0(s)\bigr) \biggr )& \leq e^{2p'^2t ||b_0||^2_{\infty}},
\label{bound2ter}
\end{split}
\end{equation}
which are bounded uniformly for $t \in [0,1]$ too. 

\subsection{Examples 2.2 and 2.3: \\temporal -and spatial- interaction}
We want now to do explicit computations for the  long-memory example
$b_i(t,\omega) = \int_0^t \epsilon(s)(\omega_i(s) - \omega_i(0))ds$ with $\epsilon$ satisfying $\eqref{contepsilon}$. The
requirement (A2) holds since
\begin{equation}
\begin{split}
& \biggl | \int_0^t \epsilon(s)(\omega_0(s)-\omega_0(0))ds - \int_0^t \epsilon(s)(\omega'_0(s)-\omega'_0(0))ds \biggr | \\
& \quad \leq 2\int_0^t \epsilon(s)ds \, \sup_{0\leq s \leq t}|\omega_0(s) - \omega'_0(s)|
\end{split}
\end{equation}
To prove the condition (A3) we first analyse the stochastic integral term
$J_t(X):= \int_0^t b_0(s,X) dX_0(s)$.
\begin{equation}
\begin{split}
J_t(X)& = \int_0^t \int_0^s \epsilon(r)(X_0(r)-X_0(0))\, dr \, dX_0(s) \\
& = \int_0^t \epsilon(r) \int_r^t \, dX_0(s)(X_0(r)- X_0(0)) \,dr \\
& = \int_0^t \epsilon(r)(X_0(t)-X_0(r))(X_0(r)-X_0(0))\, dr
\end{split}
\end{equation}
(for the interchange  of the order of integration, see for example the lecture notes \cite{Wal84}).
The integral is now an ordinary Riemann-Stieltjes one. Hence, its time-reversal satisfies
\begin{equation}
\begin{split}
J_t \circ \theta_t \,(X)  \\
& = \int_0^t \epsilon(r)(X_0(0)-X_0(t-r))(X_0(t-r)-X_0(t)) \,dr\\
& = \int_0^t \epsilon(t-r')(X_0(t)-X_0(r'))(X_0(r')-X_0(0))\, dr'.
\end{split}
\end{equation}
Similar computations lead us to the time-reversal of the functional
$$
X \mapsto  \biggl ( \int_0^t \biggl ( \int_0^s \epsilon(r)(X_0(r)-X_0(0))dr \biggr )^2 ds \biggr ).
$$
One obtains 
$$
 \int_0^t \biggl ( \int_0^s \epsilon(t-r)(X_0(r)-X_0(0))dr \biggr )^2 ds.
 $$

Thus
\begin{equation}
\begin{split}
 F^t_{0} \circ \theta_t (X)= 
& \int_0^t \epsilon(t-s)(X_0(t)-X_0(s))(X_0(s)-X_0(0))\,ds \\
& - \frac{1}{2}\int_0^t \biggl ( \int_0^s \epsilon(t-r)(X_0(r)-X_0(0))dr \biggr )^2 ds .
\end{split}
\end{equation}
As in the above example \ref{ex:Markov} the convergence of
$\exp \bigl| F^t_{0}\circ \theta_t (X)  \bigr|  $
in $\mathbb{L}_{2p}(P^x)$ for $t \rightarrow 0$ is a direct consequence of 
a uniform bound for $t \in [0,1]$ in $L^{2p'}, p'>p$, which we now prove. 
\begin{equation}
\begin{split}
 \E_{P^x} \biggl (\exp \bigl ( 2 p' \text{ } | F^t_{0} \circ \theta_t (X)| \bigr ) \biggr
)  
& \leq \E_{P^x} \biggl (\exp \bigl ( p' \varepsilon(t)(2+t \varepsilon(t))\sup_{s \leq t} \bigl [ X_0(s)-X_0(0)\bigr ]^2 \bigr ) \biggr
)\\
& \leq \E_{P^{x_0}_0} \biggl ( \sup_{s \leq t} \biggl [ \exp  \bigl( X(s)-X(0)\bigr)^2  \biggr ]^{ p' \varepsilon(t)(2+t \varepsilon(t))}\biggr
) \\
& \leq \E_{P^{x_0}_0} \biggl ( \sup_{s \leq t} \bigl [ \exp  \bigl( X(s)-X(0)\bigr)^2  \bigr ]^{ c}\biggr
) 
\end{split}
\end{equation}
where $\varepsilon (t)=: \int_0^t \epsilon(s)ds$ and $c=p' \varepsilon(1)(2+ \varepsilon(1))$.
Since $X(t)$ is a Brownian motion w.r.t. $P^{x_0}_0$, we can apply Doob's inequality 
and then obtain  
\begin{equation} \label{eq1}
\begin{split}
 \E_{P^x} \biggl (\exp \biggl ( 2 p' \text{ } | F^t_{0} \circ \theta_t (X)| \biggr ) \biggr
)  &   \leq \biggl ( \frac{c}{c-1}\biggr )^{c} \E_{P_0^{x_0}} \biggl( \exp 2c |X(t) - X(0)| \biggr) \\
& \leq \bigl( \frac{c}{c-1} \bigr)^{c} \, \E \bigl( \exp ( 2c \sqrt{t} |Z|)  \bigr)\\
& \leq \bigl( \frac{c}{c-1} \bigr)^{c} \, \E \bigl( \exp ( 2c |Z|)  \bigr) <+\infty
\end{split}
\end{equation}
where $Z$ is a standard Gaussian variable. 
The proof of (A3) is now completed for the example 2.2.\\
Example 2.3 can be treated in a very similar way, we leave the straightforward details here to the reader.

\end{document}